\documentclass[12pt,a4paper, leqno]{amsart}
\usepackage{amssymb}
\usepackage{amsmath}
\newtheorem{lemma}{Lemma}
\newtheorem{thm}{Theorem}
\theoremstyle{remark}
\newtheorem*{rem}{Remark}

\def\RR{\mathbb{R}}
\def\NN{\mathbb{N}}
\def\CC{\mathbb{C}}
\def\bl{\begin{lemma}}
\def\el{\end{lemma}}
\def\dsp{\displaystyle}
\def\ann{\operatorname{ann}}
\begin{document}
\title[Wandering domains]{An entire function with simply and
multiply connected wandering domains}
\dedicatory{Dedicated to Professor Frederick W.\ Gehring on the
 occasion of his 80th birthday}
\author{Walter Bergweiler}\thanks{Supported by the G.I.F.,
the German--Israeli Foundation for Scientific Research and
Development, Grant G-809-234.6/2003} 
\address{Mathematisches Seminar,
Christian--Albrechts--Universit\"at zu Kiel,
Lude\-wig--Meyn--Str.~4,
D--24098 Kiel,
Germany}
\email{bergweiler@math.uni-kiel.de}
\date{}
\begin{abstract}
We modify a construction of Kisaka and Shishikura to show 
that there exists an entire function $f$ which has both a simply
connected and a multiply connected wandering domain. Moreover, 
these domains are contained
in the set $A(f)$ consisting of the points where the 
iterates of $f$ tend to infinity fast.
The results answer questions by Rippon and Stallard.
\end{abstract}
\maketitle
\section{Introduction and results} \label{intro}
Let $f$ be an entire or rational function. The
{\em Fatou set}
$F(f)$ is defined as the set where the iterates $f^n$ of
$f$ form a normal family. If $U_0$ is a component of $F(f),$ then
$f^n(U_0)$ is contained in a component $U_n$ of $F(f)$. If all
$U_n$ are different, then $U_0$ is called a 
{\em wandering domain} of~$f$.
While a famous theorem of Sullivan~\cite{Sul85} says that rational
functions do not have wandering domains, it had been shown already
earlier by Baker~\cite{Bak76} that such domains may exist for transcendental
entire functions. While the wandering domain in Baker's example
was multiply connected, examples of simply connected wandering
domains were given later by various authors;
see~\cite[p.\ 106]{Her84}, \cite[p.\ 414]{Sul85},
\cite[p.\ 564, p.\ 567]{Bak84} and~\cite[p.\ 222]{Dev86}. 
Baker~\cite[Theorem~2]{Bak85}
showed that his
construction can be modified to yield wandering domains of
infinite connectivity. Recently Kisaka and Shishikura~\cite{Kis06}
constructed
an example with a multiply connected wandering domain of finite
connectivity, thereby answering a question of
Baker. In fact, they showed that the connectivity may take
any preassigned value. Here we
modify the construction of Kisaka
and Shishikura to prove the following result.
\begin{thm} \label{thm1}
There exists an entire function which has both a simply connected
and a multiply connected wandering domain. 
\end{thm}
The question whether an entire function with this property exists had been raised by
Rippon and Stallard~\cite[p.~1125, Remark~3]{Rip05}.
In the same paper, Rippon and Stallard also asked a question about
the set $A(f)$ introduced in~\cite{Ber99}. This is defined by
$$A(f):=\{z:\mbox{ there exists }L\in \NN \mbox{ such that
}|f^n(z)|>M(R,f^{n-L})\mbox{ for }n>L\},$$ 
where
$M(r,f):=\max_{|z|=r} |f(z)|$ and $R> \min_{z\in
J(f)}|z|$. Roughly speaking, $A(f)$ consists of the points $z$ where
$f^n(z)$ tends to infinity ``as fast as possible.'' Rippon and
Stallard showed that $A(f)$ has no bounded
components and that the closure of every multiply connected
wandering domain is contained in $A(f)$. They also showed that if a
simply connected wandering domain intersects $A(f),$ then it must
lie entirely in $A(f),$ and they ask~\cite[p.~1126, Remark~4]{Rip05}
whether an entire function
$f$ with such a simply connected wandering domain exists. It turns
out that an example with this property is provided by the function
constructed in Theorem~\ref{thm1}.
\begin{thm}  \label{thm2}
There exists an entire function $f$ for which $A(f)$ contains
a simply connected wandering domain.
\end{thm}
As mentioned, our
construction is largely based on that of Kisaka and Shishikura. We
state two of their lemmas in \S\ref{lemmas} and then repeat their
construction in \S\ref{construction}.
There is only one minor change in the construction, which will 
be explained at the beginning of \S\ref{proof1}.
In  the remainder of \S\ref{proof1} we 
then show that the function constructed has a simply connected
wandering domain, thereby proving 
Theorem~\ref{thm1}.
In \S\ref{proof2} we prove Theorem~\ref{thm2}.

\section{Two lemmas of Kisaka and Shishikura} \label{lemmas}
Kisaka and Shishikura first construct a quasiregular map $g:\CC
\to \CC$ and then obtain the entire function $f$ with the
following lemma.

\bl
\label{KS31}
\cite[Theorem 3.1]{Kis06} Let $g$ be a
quasiregular mapping from $\CC$ to $\CC$. Suppose that there are
(disjoint) measurable sets $E_j\subset \CC$ $ (j=1,2,...)$
satisfying: \begin{itemize}
    \item[(a)] For almost every $z\in \CC,$ the $g$-orbit of $z$ passes
    $E_j$ at most once for every~$j$;
    \item[(b)] $g$ is $K_j$-quasiregular on $E_j$;
    \item[(c)] $K_\infty :=\prod^\infty _{j=1} K_j<\infty$;
    \item[(d)] $g$ is holomorphic a.e. outside $\bigcup^\infty_{j=1}
    E_j$ (i.e. $\frac{\partial g}{\partial \bar{z}}=0$ a.e. on $\CC
    \backslash \bigcup^\infty_{j=1}E_j$).
\end{itemize}
Then there exists a $K_\infty$-quasiconformal 
map $\varphi$ such that $f=\varphi\circ g\circ \varphi^{-1}$ 
is an entire function\el
In order to construct $g$ they need to ``interpolate'' two
polynomials given on circles by a quasiregular map with small
dilatation. This is done with the following result,
where $\log$ denotes the principal branch of the logarithm.
\bl 
\label{KS63}
\cite[Lemma 6.3]{Kis06}
Let $k\in \NN$,
$b,\omega\in \CC\setminus\{0\}$ and $\rho^\sharp, \lambda^\sharp,\rho^
\flat,\lambda^\flat \in \RR$ with
$0<\lambda^\flat<\rho^\flat<1<\rho^\sharp<\lambda^\sharp$.

\begin{itemize}
    \item[(a)] Suppose that these constants satisfy
\[\rho^\sharp \geq 2|\omega|, \quad
\lambda^\sharp\geq e\rho^\sharp,\quad
C^\sharp:=1-\frac{1}{k+1}\left(\frac{|\log
b|}{\log(\lambda^\sharp/\rho^\sharp)}
+\frac{4|\omega|}{\rho^\sharp}\right)>0,\]
Then the map $bz^k(z-\omega)$
on $|z|=\rho^\sharp$ and $z^{k+1}$ on $|z|=\lambda^\sharp$ can be
interpolated on $\rho^\sharp\leq|z|\leq \lambda^\sharp$ with a
$K$-quasiregular map $g$ where $K\leq 1/C^\sharp$.
    \item[(b)] Suppose that these constants satisfy
\[|\omega|\geq 2\rho^\flat, \quad \rho^\flat\geq e\lambda^\flat, \quad 
   C^\flat:=1-\frac{1}{k}\left(\frac{|\log
(-b\omega)|}{\log(\rho^\flat/\lambda^\flat)}+\frac{4\rho^\flat}{|\omega|}\right)>0.\]Then
the map $bz^k(z-\omega)$ on $|z|=\rho^\flat$ and $z^k$ on
$|z|=\lambda^\flat$ can be interpolated on $\lambda^\flat
\leq|z|\leq \rho^\flat$ with a 
$K$-quasiregular map $g$ where $K\leq 1/C^\flat$.
\end{itemize}
\el

\section{Construction of $f$} \label{construction}
As mentioned, we follow closely the ideas of Kisaka and Shishikura
and will first construct a quasiregular map $g:\CC \to\CC$
and then obtain $f$ via Lemma~\ref{KS31}.

We denote by $\ann(r,R)$ the open annulus with inner radius $r$ and
outer radius~$R$; that is, $\ann(r,R):=\{z\in \CC:r<|z|<R\}$. The
idea is to choose sequences $(a_n)$ and $(R_n)$ such that the map
$z\mapsto a_nz^{n+1}$ maps $\ann(R_n, R_{n+1})$ onto $\ann(R_{n+1},
R_{n+2})$. The map $g$ will then be defined by $g(z)=a_nz^{n+1}$
on a large subannulus of $\ann(R_n, R_{n+1})$, and will interpolate
the mappings $z\mapsto a_{n-1}z^n$ and $z\mapsto a_nz^{n+1}$ in an
annulus containing the circle $\{z:|z|=R_n\}$.

Choosing $R_1>R_0:=1$ we obtain sequences $(R_n)$ and $(a_n)$ as
required by putting
\[R_{n+1}:=\frac{R_n^{n+1}}{R_{n-1}^n}\] and
\[a_n:= \frac{R_{n+1}}{R_n^{n+1}}=\frac{1}{R_{n-1}^{n}}.\]
Various estimates in the sequel will require that $R_1$ has been
chosen large enough. Note that with $\gamma:=\log R_1$ we have
\[\log\frac{R_{n+1}}{R_n}=n \log\frac{R_n}{R_{n-1}}=...=n!
\log\frac{R_1}{R_0}=\gamma n!.\]We define sequences
$(P_n),(Q_n),(S_n)$ and $(T_n)$ by
\[\log\frac{T_n}{S_n}=\log\frac{S_n}{R_n}=\log\frac{R_n}{Q_n}=\log\frac{Q_n}{P_n}=\sqrt{\log\frac{R_{n+1}}{R_n}}=\sqrt{\gamma n!}\]Choosing
$R_1>e$ we have $\gamma>1$ and thus
$$\frac{T_n}{S_n}=\frac{S_n}{R_n}=\frac{R_n}{Q_n}=\frac{Q_n}{P_n}>e.$$
We also have
\begin{eqnarray*}
  \log\frac{P_{n+1}}{T_n} 
  & = & 
  -\log\frac{Q_{n+1}}{P_{n+1}}  
  -\log\frac{R_{n+1}}{Q_{n+1}}  
  +  \log\frac{R_{n+1}}{R_n}  
  - \log\frac{S_n}{R_n} 
  - \log\frac{T_n}{S_n} \\
    & = & -2\sqrt{\gamma(n+1)!}  +  \gamma n!  -  2\sqrt{\gamma n!} \\
   & > & 0             \\
\end{eqnarray*}
for all $n\in\NN$, provided that $R_1$ and hence $\gamma$
is large enough. Thus
\[P_n<Q_n<R_n<S_n<T_n<P_{n+1}\] 
for all $n\in\NN$. We now define
\[b_n:=-\frac{(n+1)^2}{n+2}\left(\frac{n+1}{n}\right)^n
a_n=-\frac{(n+1)^2}{n+2}\left(\frac{n+1}{n}\right)^n\frac{R_{n+1}}{R_n^{n+1}}\]
for $n\in \NN$. We also put $E_1:=\ann(P_2, Q_2)$ and
\[E_n:= \ann (S_n,T_n) \cup \ann(P_{n+1}, Q_{n+1})\] 
for $n\geq 2$.

We shall show that there exists a quasiregular map $g:\CC
\to\CC$ with the following properties:

\begin{itemize}
    \item[(i)]$g(z)=a_1z^2$ for $|z|\leq P_2$;
    \item[(ii)]$g(z)=a_nz^{n+1}$ for $T_n\leq|z|\leq P_{n+1}$ and
    $n\geq 2$;
    \item[(iii)]$g(z)=b_n(z-R_n)z^n$ for $Q_n\leq|z|\leq S_n$ and
    $n\geq 2$;
    \item[(iv)]$g$ is $K_n$-quasiregular in $E_n$ for $n\geq1,$
    with $K_n:=1+{1}/n^2$;
    \item[(v)]$g(\ann(S_n,Q_{n+1}))\subset \ann(S_{n+1}, Q_{n+2})$ for 
    $n\geq 1$.
\end{itemize}
Since $E_n\subset \ann(S_n,Q_{n+1})$ and since the annuli $\ann(S_n,
Q_{n+1})$ are pairwise disjoint it then follows that $g$ satisfies
the hypothesis of Lemma~\ref{KS31}. Thus there exists a quasiconformal map
$\varphi$ such that $f:=\varphi\circ g\circ \varphi^{-1}$ is
entire. This function $f$ then has the desired
properties.

In order to show that a map $g$ with the properties stated exists,
we simply define~$g$ by (i), (ii), (iii) in the ranges given there
and thus have defined $g$ in $\CC\backslash \bigcup_{n=1}^\infty E_n$.

To define $g$ in $\ann(P_n,Q_n),$ where $n\geq 2,$ we consider
$G(z):=g(R_nz)/R_{n+1}$.
For $|z|=\lambda^\flat:=P_n/R_n$ we then have
\[G(z)=a_{n-1}\frac{(R_nz)^n}{R_{n+1}}=
\frac{R_n}{R_{n-1}^n} \frac{R_n^{n} z^n}{ R_{n+1}}=z^n\] 
and for $|z|=\varrho^\flat:=Q_n/R_n$ we
have
$$ G(z)  =  b_n\frac{(R_n z-R_n)(R_nz)^n}{R_{n+1}}  = 
b_n\frac{R_n^{n+1}}{R_{n+1}}(z-1)z^n 
 = c_n(z-1)z^n $$
with
\[c_n:=\frac{b_nR_n^{n+1}}{R_{n+1}}=
-\frac{(n+1)^2}{n+2}\left(\frac{n+1}{n}\right)^n.\]
Now
\[\varrho^\flat=\frac{Q_n}{R_n}=\exp (-\sqrt{\gamma n!})\]
and 
\[\lambda^\flat=\frac{P_n}{R_n}=\exp (-2\sqrt{\gamma
n!})=(\varrho^\flat)^2.\]
Thus $\varrho^\flat \geq e \lambda^\flat$
since $\gamma\geq 1$ and also $2\varrho^\flat\leq 1$. By Lemma~\ref{KS63}, (b),
there exists a $K$-quasiregular map
$G_n:\{z\in \CC:\lambda^\flat\leq|z|\leq\varrho^\flat\}
\to\CC$
such that $G_n(z)=z^n$ for $|z|=\lambda^\flat$ and
$G_n(z)=c_n(z-1)z^n$ for $|z|=\varrho^\flat,$ with $K\leq 1/C^\flat$
where
\[C^\flat:=1-\frac{1}{n}\left(\frac{|\log(-c_n)|}{\log(\varrho^\flat/\lambda^\flat)}+4\varrho^\flat\right),\]
provided that $C^\flat>0$.
But since
$$|\log(-c_n)|  =  \log\left(\frac{(n+1)^2}{n+2}\left(\frac{n+1}{n}
\right)^n\right) 
 \leq  \log((n+1)  e) =  1+ \log(n+1) $$
we may in fact achieve that  
\[C^\flat\geq 
1-\frac{1}{n}\left(\frac{1+\log(n+1)}{\sqrt{\gamma n!}}+4\exp\left(
{{-\sqrt{\gamma n!}}}\right)\right)
\geq 
1-\frac{1}{(n-1)^2+1}\] 
for all $n\geq 2$ by choosing $\gamma$ large enough. Thus
\[K\leq \frac{1}{C^\flat} \leq 1+
\frac{1}{(n-1)^2}.\]
Putting
\[g(z):=R_{n+1}G_n\left(\frac{z}{R_n}\right)\] 
for $z\in \ann(P_n,Q_n)\subset E_{n-1}$ we
see that (iv) holds for $z\in E_{n-1}\cap \ann(P_{n+1},Q_{n+1})$.
Similarly we define $g$ in the remaining part of $E_n;$ that is,
in $E_n\cap \ann(S_n,T_n)$. Here we use the first part of Lemma~\ref{KS63}.

To prove (v) we note that if $z\in \ann(S_n,Q_{n+1}),$ then, by
the maximum principle,
\begin{eqnarray*}
   |g(z)|
   & \leq &  \max_{|\zeta|=Q_{n+1}}|g(\zeta)|\\
   & = & \max_{|\zeta|=Q_{n+1}}\left|b_{n+1}
         (\zeta-R_{n+1})\zeta^{n+1}\right| \\
   & = & \left|b_{n+1}\right|\left(R_{n+1}+Q_{n+1}\right)Q_{n+1}^{n+1} \\
   & = & \frac{(n+2)^2}{n+3}\left(\frac{n+2}{n+1}\right)^{n+1}
         R_{n+2} \left(1+\frac{Q_{n+1}}{R_{n+1}}\right)
        \left(\frac{Q_{n+1}}{R_{n+1}}\right)^{n+1} \\            
   & \leq & 2e(n+2)R_{n+2}\left(\frac{Q_{n+1}}{R_{n+1}}\right)^{n+1} \\
   & = & 2e(n+2) R_{n+2} 
   \exp\left({-(n+1)\sqrt{\gamma(n+1)!}}\right) \\
   & \leq & R_{n+2}\exp\left({-\sqrt{\gamma(n+2)!}}\right) \\
   & = & Q_{n+2}
\end{eqnarray*}
if $\gamma$ is large enough. Similarly, noting that
$g$ has no zeros in $\ann(S_n,Q_{n+1})$ and using the
minimum principle, we find that
$|g(z)|\geq S_{n+1}$ for $z\in \ann(S_n,Q_{n+1})$. We deduce that (v)
holds.

As in the paper of Kisaka and Shishikura we deduce from (v) that $g^n(z)
\to \infty$ as $n\to \infty$ for $z\in
\ann(S_1,Q_2),$ while $g(0)=0$ by (i). This implies that 
$\varphi(\ann(S_1,Q_2))$ lies in a multiply connected
component $U_1$ of the Fatou set of $f$,
with $f^n|_{U_1}\to\infty$ as $n\to\infty$.
Since multiply connected components of the Fatou set are always
bounded by a result of Baker~\cite[Theorem~1]{Bak75},
this implies that $f$ has a
multiply connected wandering domain.

\begin{rem}
The Fatou set and the other concepts of complex dynamics can
also be defined for quasiregular maps, by carrying over the definitions
from the holomorphic case literally.
In order to retain the basic features of the theory one
has to require, however, that all
iterates of $g$ are $K$-quasiregular with the same~$K$. 
Such maps are called  {\em uniformly quasiregular}.
That our function 
$g$ is uniformly quasiregular
follows directly from the definition of $g$, or
trivially from the representation $g=\varphi^{-1} \circ f
\circ\varphi$.
Lemma~\ref{KS31} says, essentially, that uniformly quasiregular selfmaps 
of the plane are
quasiconformally conjugated to entire functions;
see also~\cite{Gey,Hin}
for this result.
\end{rem}

\section{Proof of Theorem 1: $f$ has a simply connected wandering domain}
\label{proof1}
The sequence $(\xi_n)$ of critical points of $g$ is given by
$\xi_1:=0$ and $\xi_n:=n R_n/(n+1)$ for $n\geq 2$. The only
difference between the present construction and that of Kisaka and
Shishikura concerns the orbits of these points. While we have
chosen the values $b_n$ such that
$$ g(\xi_n)  =  b_n(\xi_n-R_n)\xi_n^n 
= -b_n \frac{R_n}{n+1}\left(\frac{n}{n+1}\right)^n R_n ^n 
=\frac{n+1}{n+2}R_{n+1} 
=\xi_{n+1}$$
for $n\geq 2$,
Kisaka and Shishikura worked with different values
of $b_n$ which yielded $g(\xi_n)=R_{n+1}$ and hence
$g^2(\xi_n)=0$.

Denote by $D(a,r)$ the disk of radius $r$ around~$a$. Let
$\delta>0$ and define $D_n:=D(\xi_n, \delta R_n/n^4)$ for $n\geq
2$. We shall show that if $\delta$ is sufficiently small, then
$g(D_n)\subset D_{n+1}$ for all~$n$. This implies that $D_n\subset
F(g)$ for $n\geq 2$. We will then show that $D_n$ lies in a simply
connected wandering domain of $g$ and thus $\varphi(D_n)$ lies in
a simply connected wandering domain $V_n$ of~$f$.
Moreover, we will see in \S\ref{proof2}
that $V_n\subset A(f)$ for all~$n$.

First we note that if $\delta$ is small enough, then $D_n\subset\ann(Q_n,R_n)$ 
so that $g$ is holomorphic in $D_n$ and
\begin{eqnarray*}
  |g''(z)| & = & \left|b_n\left(n(n+1)z^{n-1}-R_n n(n-1)z^{n-2}\right)\right|\\
   & \leq & |b_n|(n(n+1)+n(n-1))R_n^{n-1} \\
   & = & 2n^2|b_n|  R_n^{n-1} 
\end{eqnarray*}
for $z\in D_n$. Thus
\begin{eqnarray*}
  |g'(z)| & = & |g'(z)-g'(\xi_n)| \\
   & \leq & \int_{\xi_n}^z |g''(\zeta)|\;|d\zeta| \\
   & \leq & 2n^2  |b_n| R_n^{n-1} \frac{\delta R_n}{n^4} \\
   & = & \frac{2\delta}{n^2}|b_n| R_n^n 
\end{eqnarray*}
for $z\in D_n$. It follows that if $z\in D_n,$ then
\begin{eqnarray*}
  |g(z)-\xi_{n+1}| & = & |g(z)-g(\xi_n)| \\
   & \leq  & \int^z_{\xi_n} |g'(\zeta)|\;|d\zeta| \\
   & \leq &  \frac{2\delta}{n^2}|b_n| R_n^n \frac{\delta R_n}{n^4}\\
   & = & \frac{2\delta^2}{n^6}\frac{(n+1)^2}{n+2}
   \left(\frac{n+1}{n}\right)^n R_{n+1} \\
   & \leq & \frac{2\delta^2 e (n+1)}{n^6} R_{n+1}\\
   & \leq & \frac{\delta R_{n+1}}{(n+1)^4}, 
\end{eqnarray*}
provided $\delta$ is sufficiently small. Thus
$g(D_n)\subset D_{n+1}$ for $n\geq 2$. As already mentioned, this
implies that $D_n$ lies in a  Fatou component $V'_n$ of $g$ and
thus $\varphi(D_n)$ lies in the Fatou component
$V_n:=\varphi(V'_n)$ of~$f$. By $U'_n$ we denote the multiply
connected Fatou component of $g$ which contains
$\ann(S_n,Q_{n+1}),$ and by $U_n:=\varphi(U'_n)$ the corresponding
Fatou component of~$f$. As mentioned at the end of \S\ref{construction},
the $U_n$ are wandering domains. In fact, we have $U_m\neq U_{l}$ for
$m\neq l$ (and thus $U'_m\neq U'_l$ for $m\neq l$). We
shall show that $U'_m\neq V'_l$ (and thus $U_m\neq V_l$) for
all $m$ and~$l$. Since $V'_l$ lies ``between'' the annuli
$\ann(S_{l-1},Q_l)$ and $\ann(S_l,Q_{l+1})$ and thus
``between'' the domains 
$U'_{l-1}$ and $U'_l,$ it suffices to show that
$U'_m\neq V'_l$ for $l=m$ and $l=m+1$.

Suppose that $U'_m=V'_l$ where $l=m$ or $l=m+1$. Since $D_l\subset V_l'$ and
$\ann(T_m,P_{m+1})\subset U_m'$ there
exists a simply connected domain $\Omega_m$ with
\[D_l\cup \ann(T_m,P_{m+1})\backslash
(-P_{m+1},T_m)\subset\Omega_m\subset U'_m.\]
Since $g^n(z)\neq 0$
for $z\in\Omega_m$ and since the $g^n$ are $K$-quasiregular for
some $K$, we may define for $n>m$ a $K$-quasiregular map
$h_n:\Omega_m\to\CC$ by
\[h_n(z)=\left(\frac{g^{n-m}(z)}{R_n}\right)^{m!/n!},\]
for some branch
of the root. We will show that the $h_n$ form a normal family so
that $(h_n)$ has a convergent subsequence, say $h_{n_k}\to
h$. Next we will show that $h_n(z)=h(z)=z/R_m$ for $z\in
\ann(T_m,P_{m+1})\backslash (-P_{m+1},T_m)$, if the branch of
the root has been suitably chosen. This implies in particular that $h$ is nonconstant. On
the other hand, we will see that $h$ is constant in $D_l$ so
that we obtain a
contradiction.

To prove that $(h_n)$ is normal we note that if $z\in
\Omega_m\subset U'_m,$ then $g^{n-m}(z)\in U'_n$ and thus
$|g^{n-m}(z)|\leq S_{n+1},$ since $\ann(S_{n+1},Q_{n+2})\subset
U'_{n+1}$. Hence
\[
|h_n(z)|
\leq
\left|\frac{S_{n+1}}{R_n}\right|^{m!/n!}
=
\left|\frac{S_{n+1}}{R_{n+1}}\frac{R_{n+1}}{R_n}\right|^{m!/n!}
\]
for $z\in \Omega_m$. We deduce that
\begin{eqnarray*}
  \log|h_n(z)|  
  & \leq & \frac{m!}{n!}
  \left(\log \frac{S_{n+1}}{R_{n+1}}+\log \frac{R_{n+1}}{R_n}\right) \\
   & = & \frac{m!}{n!} \left(\sqrt{\gamma(n+1)!}+\gamma n!\right)\\
   & \leq & 2\gamma m! 
\end{eqnarray*}
for $z\in \Omega_m$ and large $n,$ and this yields the desired
normality.

It is not difficult to see by induction that if 
$z\in \ann(T_m,P_{m+1})\backslash(-P_{m+1},T_m)$ and $n>m$,
then 
\[
g^{n-m}(z)  = R_n\left(\dsp\frac{z}{R_m}\right)^{n!/m!} 
\in  \ann(T_n,P_{n+1})\backslash(-P_{n+1},T_n)
\]
so that
\[h_n(z)=h(z)=\frac{z}{R_m}\]
if the branch of the root in the
definition of $h_n$ has been suitable chosen.
In particular, $h$ is nonconstant.

For $z\in D_l$ we have
\[g^{n-m}(z)\in D_{n-m+l}\]
and thus $g^{n-m}(z)\in D_n$ or $g^{n-m}(z)\in D_{n+1}$ depending
on whether $l=m$ or $l=m+1$. In the first case we have
\[\frac{g^{n-m}(z)}{R_n}\in D\left(\frac{\xi_n}{R_n},
\frac{\delta}{n^4}\right)=D\left(\frac{n}{n+1},\frac{\delta}{n^4}\right)\subset
D\left(1,\frac{1}{2}\right)\] for large $n,$ and this implies that
$h$ is constant in $D_l$. In the second case, that is for $l=m+1,$
we have
\[\frac{g^{n-m}(z)}{R_n}\in
D\left(\frac{\xi_{n+1}}{R_n},\frac{\delta
R_{n+1}}{(n+1)^4R_n}\right)\] and hence
\[\left(1-\frac{2}{n}\right)\frac{R_{n+1}}{R_n}\leq 
\left|\frac{g^{n-m}(z)}{R_n}\right|\leq \frac{R_{n+1}}{R_n}\]
for large~$n$.
Since $\sqrt[n!]{1-2/n}\to 1$ as $n\to\infty$ and
\[\left|\frac{R_{n+1}}{R_n}\right|^{m!/n!}
=\exp (\gamma m!)\] 
we deduce
that $|h(z)|=\exp (\gamma m!)$ for $z\in D_l$. Thus $h$ is
constant in $D_l$ in this case as well. As already noted, this
is a contradiction. This completes the proof that $U'_m\neq
V'_l$ and hence $U_m\neq V_l$ for all $l$ and~$m$.

It remains to show that $V_l$ is simply connected for all
$l$. Suppose that this is not the case. It is not difficult to
show that then all $V_l$ are multiply connected. By a result
of Baker~\cite[Theorem~3.1]{Bak84} 
there exists $k\in \NN$ and a Jordan curve $\tau$ in
$V_k$ whose interior contains 0. This implies that $D(0,Q_k)$ is
contained in the interior of $\tau':=\varphi^{-1}(\tau)$. The
argument principle implies that the winding number of $g(\tau')$
around 0
is at least~$k$. 
This winding number is equal to
the winding number of $f(\tau)$ around 0, 
and thus the latter winding number is also at least~$k$.
Induction shows that the winding number of
$\tau_n:=f^{n-k}(\tau)$ around 0 is at least $(n-1)!/(k-1)!$
for $n\geq k$. A contradiction will now be obtained from a 
consideration of the hyperbolic length of
$\tau_n$ in $V_n$. We denote the hyperbolic length of a curve
$\sigma$ in a hyperbolic domain $U$ by $\ell(\sigma,U)$. By the
Schwarz-Pick-Lemma we have
\[\ell(\tau_n,V_n)\leq \ell(\tau,V_k)\] for all $n\geq k$. On
the other hand, we have $V_n'\subset \ann(Q_n,S_n)$ and thus $V_n
\subset \varphi(\ann(Q_n,S_n))$. This implies that
\[\ell(\tau_n,V_n)\geq \ell(\tau_n,\varphi(\ann(Q_n,S_n))).\]
Now $\ann(Q_n,S_n)$ is an annulus of modulus $\log(S_n/Q_n)/(2\pi)$. 
Since $\varphi$ is $K$-quasi\-con\-for\-mal this yields
that $\varphi(\ann(Q_n,S_n))$ has modulus at most
\[\dsp\frac{K}{2\pi}\log
\left(\frac{S_n}{Q_n}\right)=
\frac{1}{2\pi}\log\left(\left(\frac{S_n}{Q_n}\right)^K\right).\]
It follows that there exists a conformal map
$\psi:\varphi(\ann(Q_n,S_n)) \to \ann(1,r_n)$ where $r_n\leq
(S_n/Q_n)^K$. 
We may choose $\psi$ such that $|\psi(\varphi(z))|\to 1$ as 
$|z|\to Q_n$.
Put $\sigma_n:=\psi(\tau_n)$. Then
\[\ell(\sigma_n,
\ann(1,r_n))=\ell(\tau_n,\varphi(\ann(Q_n,S_n)))\] and the winding
number of $\sigma_n$ around $0$ is the same as that of $\tau_n$ and
thus at least $(n-1)!/(k-1)!$. 
We note  that the density $\varrho(z)$ of
the hyperbolic metric in $\ann(1,r_n)$ is given by 
(see, e.~g.,~\cite[p.~12]{McM94})
\[\varrho(z)=\dsp\frac{\pi}{|z| \sin(\pi\log |z|/\log r_n) \log
r_n}.\] 
In particular we have $\varrho(z)\geq \pi/(|z|\log r_n)$ and thus we conclude that
$$ \ell(\sigma_n, \ann(1,r_n))  =  \int_{\sigma_n}\varrho(w)|dw| 
 \geq  \frac{\pi}{\log r_n}\int_{\sigma_n} \frac{|dw|}{|w|}
\geq  \frac{2\pi^2}{\log r_n} \frac{(n-1)!}{(k-1)!}. $$ 
Since $\log r_n\leq K \log \left({S_n}/{Q_n}\right)
=2K\sqrt{\gamma n!}$ we
deduce that
\[\ell(\sigma_n, \ann(1,r_n))\geq 
\frac{\pi^2}{K(k-1)!}\frac{(n-1)!}{\sqrt{\gamma
n!}}\] so that
\[\ell(\sigma_n, \ann(1,r_n))\to \infty\] as $n\to
\infty$.

On the other hand, our previous estimates imply that
\[\ell(\sigma_n, \ann(1,r_n))=\ell(\tau_n,\varphi(\ann(Q_n,S_n)))
\leq \ell(\tau_n,V_n)\leq
\ell(\tau,V_k).\] This is a contradiction. Thus $V_\ell$ is
simply connected for all $\ell$. This completes the proof of
Theorem~\ref{thm1}.

\begin{rem}
Except for the fixed point $\varphi(0)$,
the critical points of $f$ are contained 
in the simply connected wandering domains $V_n$ of $f$. 
Thus the wandering domains $U_n$ do not contain critical
points. Using this it can be shown with the arguments of 
Kisaka and Shishikura (in particular,~\cite[Proposition~4.5]{Kis06}) 
that the $U_n$ are doubly connected.

\end{rem}
\section{Proof of Theorem 2: the $V_k$ are in $A(f)$}
\label{proof2}
We will use the following characterisation of the set $A(f)$ given
by Rippon and Stallard~\cite[Lemma~2.4]{Rip05}. 
Here we denote for a domain $U$ by
$\tilde{U}$ the union of $U$ and its bounded complementary
components.

\bl 
\label{ripsta}
Let $f$ be a transcendental entire function and let $D$ be a
domain intersecting the Julia set of~$f$. Then
\[A(f)=\{z:\mbox{ there exists }L\in\NN \mbox{ such that }
f^{n+L}(z)\not\in \widetilde{f^n(D)}\mbox{ for }n\in \NN\}.\]\el
We apply this result with $D:=U_1$. It follows from the maximum
principle that $\widetilde{f^n(U_1)}=\widetilde{U_{n+1}}$. Since 
$U_m'\subset D(0,S_{m+1})$ we have 
$V_{l}'\cap \widetilde{U_{m}'}=\emptyset$ and hence 
$V_{l}\cap \widetilde{U_{m}}=\emptyset$ for $l\geq m+2$.
Thus we see that
if $k\geq 2$ and $z\in V_k$ so that  $f^{n+1}(z)\in V_{k+n+1}$, then 
\[f^{n+1}(z)\not\in \widetilde{U_{n+1}}=\widetilde{f^n(U_1)}.\]
Choosing $L:=1$ in Lemma~\ref{ripsta} we see
that $V_k\subset A(f)$ for $k\geq 2$. This completes the proof of
Theorem~\ref{thm2}.

\end{document}